\documentclass[11pt]{article}

\usepackage{enumerate}
\usepackage[OT1]{fontenc}
\usepackage{amsmath,amssymb}
\usepackage{natbib}
\usepackage[usenames]{color}
\usepackage{fullpage}
\usepackage{multirow}

\usepackage[colorlinks,
            linkcolor=red,
            anchorcolor=blue,
            citecolor=blue
            ]{hyperref}

\def\tr{\mathop{\text{tr}}\kern.2ex}

\long\def\comment#1{}

\def\tr{\mathop{\text{Tr}}}

\def\cS{{\mathcal{S}}}

\newcommand{\bel}{\begin{eqnarray}\label}
\newcommand{\eel}{\end{eqnarray}}
\newcommand{\bes}{\begin{eqnarray*}}
\newcommand{\ees}{\end{eqnarray*}}

\def\tr{\mathop{\text{tr}}\kern.2ex}

\long\def\comment#1{}

\def\tr{\mathop{\text{Tr}}}

\def\cS{{\mathcal{S}}}

\def\cN{{\mathcal{N}}}

\def\cB{{\mathcal{B}}}
\def\tr{{\text{Tr}}}

\usepackage{chngcntr}
\usepackage{apptools}
\usepackage{mathrsfs}

\usepackage{enumitem}
\setlist[itemize]{leftmargin=2em}

\usepackage{mathrsfs}

\def\sep{\;\big|\;}

\usepackage{breakcites}

\usepackage{chngcntr}
\usepackage{apptools}

\usepackage{microtype}
\def\##1\#{\begin{align}#1\end{align}}
\def\$#1\${\begin{align*}#1\end{align*}}
\usepackage{enumitem}

\usepackage{mathrsfs}

\def\sep{\;\big|\;}

\usepackage{enumitem}

\newcommand{\vertiii}[1]{{\left\vert\kern-0.25ex\left\vert\kern-0.25ex\left\vert #1 
    \right\vert\kern-0.25ex\right\vert\kern-0.25ex\right\vert}}
\usepackage{hyperref}
\usepackage{cleveref}
\crefname{section}{§}{§§}
\Crefname{section}{§}{§§}

\usepackage{chngcntr}
\usepackage{apptools}

\usepackage[misc]{ifsym}
\usepackage{amsmath}
\usepackage{amsthm}
\usepackage{amssymb}
\usepackage{amsopn}
\usepackage{graphicx}
\usepackage{color}
\usepackage{booktabs}
\usepackage{multicol}
\usepackage{multirow}
\usepackage{caption}
\usepackage{appendix} 
\usepackage{algorithm}
\usepackage{algorithmic}

\usepackage{epstopdf}
\usepackage{tikz}

\usepackage{bm}

\usepackage{listings}
\usepackage{matlab-prettifier}
\usepackage{float}

\usepackage{authblk}
\usepackage{parskip}

\def\rmin{\mathop{\rm minimize}}
\def\rst{\mathop{\rm subject\ to}}

\def\cN{{\cal N}}

\def\cB{{\cal B}}

\def\cS{{\cal S}}
\def\cW{{\cal W}}
\def\cV{{\cal V}}
\def\cA{{\cal A}}
\def\cN{{\cal N}}
\def\cR{{\cal R}}

\def\bx{{\bar x}}

\def\Re{{\rm I\!R}}

\def\bw{{\bar w}}
\def\bx{{\bar x}}
\def\by{{\bar y}}

\usepackage{caption}
\usepackage{subcaption}

\usepackage{amsmath}

\newtheorem{thm}{Theorem}%
\newtheorem{prop}{Proposition}%
\newtheorem{defi}{Definition}%
\newtheorem{alg}{Algorithm}%
\newtheorem{assu}{Assumption}%

\def\##1\#{\begin{align}#1\end{align}}
\def\$#1\${\begin{align*}#1\end{align*}}
\usepackage{enumitem}

\newcommand\extrafootertext[1]{%
    \bgroup
    \renewcommand\thefootnote{\fnsymbol{footnote}}%
    \renewcommand\thempfootnote{\fnsymbol{mpfootnote}}%
    \footnotetext[0]{#1}%
    \egroup
}

\makeatletter
\def\@fnsymbol#1{\ensuremath{\ifcase#1\or 1 \or 2\or
   3\or 4\or 5\or 6\or 7
   \or 8 \else\@ctrerr\fi}}
    \makeatother

\title{{\textbf{\Large{A Cardinality-Constrained Approach to Combinatorial
Bilevel Congestion Pricing}}}}

\author
{
\normalsize
Lei Guo$^{1,\dagger}$ 
\qquad Jiayang Li$^{2,\dagger}$
\qquad Yu (Marco) Nie$^{3,\dagger, \star}$
\qquad Jun Xie$^{4, \dagger}$
}

\date{}

\begin{document}

\extrafootertext{$^\dagger$Alphabetical order. $^\star$Corresponding author. Email: \texttt{y-nie@northwestern.edu}. $^1$School of Business, East China University of Science and Technology. $^2$Department of Data and Systems Engineering, The University of Hong Kong. $^3$Department of Civil and Environmental Engineering, Northwestern University. $^4$School of Transportation and Logistics, Southwest Jiaotong University.
}

\maketitle

\begin{abstract}
{
Combinatorial bilevel congestion pricing (CBCP), a variant of the mixed (continuous/discrete) network design problems, seeks to minimize the total travel time experienced by all travelers in a road network, by strategically selecting toll locations and determining toll charges.   Conventional wisdom suggests that these problems are intractable since they have to be formulated and solved with a significant number of integer variables. Here, we devise a scalable local algorithm for the CBCP problem that guarantees convergence to an approximate Karush-Kuhn-Tucker point.  Our approach is novel in that it eliminates the use of integer variables altogether, instead introducing a cardinality constraint that limits the number of toll locations to a user-specified upper bound.  The resulting bilevel program with the cardinality constraint is then transformed into a block-separable, single-level optimization problem that can be solved efficiently after penalization and decomposition.  We are able to apply the algorithm to solve, in about 20 minutes, a CBCP instance with up to 3,000 links.  To the best of our knowledge, no existing algorithm can solve CBCP problems at such a scale while providing any assurance of convergence. 
}

\end{abstract}

\section{Introduction}\label{intro} %

For nearly a century, congestion pricing has been hailed as an effective approach to managing traffic congestion in big cities \citep{pigou1920economics,vickrey1969congestion}. No event better illustrates its promises and controversies than the flip-flop of the New York City on its highly anticipated congestion pricing program \citep{economist2024}. The economic theory behind congestion pricing is sound and intuitive: {traffic congestion creates an externality, which pricing can internalize by making travelers bear the marginal, rather than average, cost of their trip} \citep[e.g.,][]{arnott1994,ferrari1995,bergendorff1997,hearn1998solving,yangandhuang2005, lindsey2006}.  Since tolls cannot be charged on all roads in a network, what one aims to achieve in practice is usually a second-best policy, which, unlike a first-best policy, may only reduce but not eliminate the externality-induced efficiency losses \citep[e.g.,][]{verhoef2002,yangandhuang2005}.

A typical second-best congestion pricing problem sets tolls on a prespecified set of roads according to a system objective. Travelers are assumed to be self-interested in that they always try to minimize their combined cost of travel time and toll when choosing the route of their trips.  %
The second-best congestion pricing problem is a leader-follower, or Stackelberg, game in which the leader is a toll-setting authority and the followers are travelers who can be seen as playing a congestion game. Such a Stackelberg congestion game belongs to a broader class of bilevel problems \citep{migdalas1995bilevel,yang1997,patriksson2002} known to be NP-hard \citep{ben1990}  --- the upper-level problem sets tolls, according to which the lower-level determines travelers' individually optimal choices (often referred to as their best response). %
Accordingly, the second-best congestion pricing problem is hereafter referred to as the bilevel congestion pricing (BCP) problem.

The BCP problem is intrinsically non-convex due to the hierarchical structure \citep[e.g.,][]{dempereview}.  It is also defined over a graph representing a road network, which further increases complexity, especially when the size of the network is large. Note that real networks used for planning purposes could easily have tens of thousands of road segments (links) and millions of origin-destination (O-D) pairs. As a result, solving the BCP problem for real-world planning or design exercises has been a long-standing computational challenge \citep{li2022differentiable}, and until recent, {no real breakthrough} has been reported despite numerous attempts  \citep[e.g.,][]{yan1996,labbe1998,yin2000genetic,lawphongpanich2004,yangandhuang2005,de2011,fallah2018,Guo2024-joc}.

The BCP problem takes the set of toll links as given and focuses on the optimization of toll levels on those links.   When the set itself needs to be optimized, we have a joint optimization problem that is in essence combinatorial. Solving this joint problem is highly desirable, as it provides a rigorous framework for guiding the selection of toll links.  In practice, a planner may first eliminate links deemed unsuitable for tolling and then rely on the joint optimization to determine which of the remaining ``feasible" links should be tolled. This approach is particularly valuable since  the number of toll links is often constrained by a budget for constructing toll facilities.  However, despite its practical appeal, this \textit{combinatorial bilevel congestion pricing} (CBCP) problem presents an even more formidable computational challenge.  If the number of feasible links is $m$ and the budget dictates the number of toll links cannot exceed $t$,  for example, there are in total $m!/(t!(m-t)!)$ ways to form the set. A brute-force approach to finding the optimal set with $t$ links would require solving $m!/(t!(m-t)!)$ BCP problems, which, as mentioned earlier, are themselves NP-hard.    Due to this immense complexity, few serious attempts have been made to solve the problem \textit{optimally}. %
Existing methods rely on either meta-heuristics without any convergence guarantee \citep{verhoef2002,yang2003,shepherd2004,harks2015computing} or approximations that provide only loose lower bounds \citep{ekstrom2012}.  The present study is motivated by this challenge.  Specifically, we aim to develop a convergent algorithm capable of locating a stationary point for the CBCP problem.

By employing a cardinality function, rather than binary variables, to code the toll state of all links, our approach avoids dealing with a mixed-integer nonlinear bilevel formulation, hence the intractability that often comes with it. The cardinality function maps a vector to the number of its non-zero components, thus transforming a combinatorial problem (i.e., picking at most $t$ toll links out of $m$) into a simple cardinality constraint (i.e., the value of the cardinality function must not exceed $t$).   This idea has been widely used in portfolio selection \citep{bienstock1996,bertsimas2009,gao2013,zheng2014}, though less known in transportation.

Our approach thus transforms the CBCP problem into a bilevel program with a special cardinality constraint.  However, the cardinality constraint, being neither convex nor continuous, is a nasty constraint, and to the best of our knowledge, no existing research has attempted to handle such a constraint in the context of bilevel programming. To overcome this hurdle, we rely on two key insights: (i) the projection onto the cardinality constraint can be computed in closed form; and (ii) the lower-level problem can be turned into an upper-level constraint defined by a gap function, which is block-wise convex with respect to the upper-level and the lower-level decision variables.  Leveraging these properties, we transform the CBCP problem into a block-separable, single-level optimization problem that can be decomposed into two simpler subproblems, for which efficient algorithms exist. Under mild conditions, we prove that the proposed algorithm guarantees convergence to an approximate Karush-Kuhn-Tucker (KKT) point of the CBCP problem. Our numerical experiments demonstrate that the algorithm can handle CBCP problems of a scale previously thought unmanageable, especially when a certain level of quality assurance is required.

In a nutshell, our study contributes to the rich literature on bilevel programming by 
\begin{itemize}
    \item[1.] proposing a novel approach to formulating a combinatorial bilevel problem, centering on a cardinality constraint that limits the number of links on which tolls can be leveled; %
    
    \item[2.] designing an efficient  algorithm that takes advantage of block-separability of the new formulation, which, unlike the heuristics that have been the mainstay for solving the CBCP problem,  offers a guarantee of convergence toward a KKT point;

    \item[3.] demonstrating the computational superiority of the proposed algorithm by applying it to large-scale CBCP problems and comparing it with benchmark heuristics.

\end{itemize}

The remainder of the paper is organized as follows. In Section \ref{sec:liter}, we briefly review related studies, and in Section \ref{setting}, we present the problem setting. 
Section \ref{model} introduces the newly proposed CBCP model, and Section \ref{sec:PD} presents the decomposition method by exploiting the revealed structure. Section \ref{sec:num} evaluates the newly proposed model and algorithm numerically through extensive experiments conducted on publicly available test networks. Section \ref{sec:con} presents conclusions and future work. All proofs are included in the appendix.

\section{Related studies}\label{sec:liter}

The combinatorial bilevel congestion pricing (CBCP) problem that chiefly concerns us here may be viewed as an extension of the BCP problem widely studied in the literature on bilevel programming.  Hence, we shall begin with that literature, before turning to the studies that specifically address the BCP problem, as well as its combinatorial variant.

\subsection{Bilevel programming}

A bilevel program is a hierarchical optimization problem whose feasible region is partly determined by the solution to another optimization problem \citep{bracken1973}. Its history can be traced back to \cite{Stackelberg}, who formulated the asymmetric competition between two players as a leader-follower game.  Bilevel programming has seen broad applications, since the hierarchy underlies many complex socio-technical systems across many domains, from economics and planning to engineering and computer science \cite[e.g.,][]{dempereview}. Despite their popularity and significant potential, bilevel programs are notoriously difficult to solve, primarily due to the inherent non-convexity introduced by their hierarchical structure \citep{ben1990}.
The past thirty years have seen considerable efforts dedicated to solving bilevel programs efficiently; see \cite{Vicente-Calamai,dempereview,colson} for reviews of these efforts. 

In a linear bilevel program, where both the upper-level and lower-level problems can be formulated as linear programs, the optimal solution occurs at an extreme point of the constraint set. This property allows for the development of exact algorithms, such as the branch-and-bound algorithm \citep{tuy1998}  and the k\textit{th}-best algorithm \citep{bialas1984}. For a general bilevel program, the standard approach is to transform it to a single-level problem  by reformulating the lower-level problem \citep{dempereview}. When the lower-level problem's best response function is single-valued and continuously differentiable, substituting it into the upper-level objective function yields an implicitly defined composite function that is also continuously differentiable. The problem can then be solved by  gradient descent methods provided the gradient can be computed efficiently  \citep[e.g.,][]{kolstad1990,savard1994,falk1995}. This approach is particularly popular in transportation applications, see \citet{li2022differentiable} and \citet{li2023achieving} for some recent efforts attempting to improve its scalability. It has also seen applications in machine learning of late \cite[see, e.g.,][for an overview]{liu2021investigating}. When the lower-level problem is a convex program, the bilevel program can be reformulated as a mathematical program with equilibrium constraints (MPEC), which has been extensively studied in the literature. Many algorithms have been proposed for MPECs, ranging from the regularization method \citep[e.g.,][]{hoheisel2013theoretical} and the penalty function method \citep[e.g.,][]{hu2004convergence} to the smoothing method  \citep[e.g.,][]{facchinei1999} and the augmented Lagrangian method \citep[e.g.,][]{izmailov2012}. Recently, \cite{guo2021} and \cite{guo2024} combined regularization and smoothing techniques to solve MPECs with a non-Lipschitz upper-level objective function.  A third approach is to reformulate the bilevel program as a single-level nonlinear program by employing the value function of the lower-level problem, an idea that was first proposed by \cite{outrata1990} and then exploited by a few other authors \cite[e.g.,][]{ye1995,meng2001equivalent}. {The value function may not be continuously differentiable. When this is the case, it needs to be smoothed, leading to the so-called smoothing approximation method \citep[e.g.,][]{lin2014}.  Otherwise, gradient descent methods can be directly applied \cite[e.g.,][]{meng2001equivalent,liu2022}.}

Few studies have considered cardinality constraints in bilevel programs. The only exception we are aware of is \cite{aussel2024cardinality}. Their proposal was to reformulate the constraint as a complementarity condition, and combine it with similar conditions derived from the lower-level problem to form an MPEC, which is then solved using standard techniques. Thus, \cite{aussel2024cardinality} made no attempt to contribute to algorithm development. Also worth noting is that their approach may not be applicable in our context—certainly not for large-scale instances—since directly solving an MPEC with a large number of complementarity constraints is computationally intractable

\subsection{Bilevel congestion pricing}
The majority of the literature on the BCP problem assumed the set of toll links is already given, hence leaving out the combinatorial optimization part of the problem.  \citet{yan1996} proposed to solve this version of the BCP problem using an implicit function-based formulation, for which the gradient of the composite function is computed using a sensitivity analysis-based method. {In \citet[Chapter 5]{lim2002}, an MPEC approach is proposed to tackle the BCP problem. To address the computational challenges posed by a large number of complementarity constraints, an iterative procedure is developed to progressively generate paths, corresponding to the vertices of a polyhedron that forms the feasible region of the lower-level traffic assignment problem.} \citet{lawphongpanich2004} turned the lower-level problem into a variational inequality constraint in the upper level, and solved the resulting single-level problem by devising a cutting plane algorithm. In \citet[Chapter 5][]{yangandhuang2005}, a value function-based reformulation is proposed and solved by a standard augmented Lagrangian algorithm. Some have also attempted to solve the BCP problem using heuristics that offers no convergence assurance. One example is \citet{ferrari2002}, who proposed a downhill simplex algorithm. Another is the genetic algorithm \citep{yin2000genetic}. \citet{fallah2018} compared the performance of the genetic algorithm and another meta-heuristic called particle swarm optimization in solving the BCP problem.   

While scarce, there have been studies attempting to solve the combinatorial BCP (CBCP) problem. \cite{verhoef2002} discussed a few heuristic strategies that hinge on ranking links according to certain selection criteria. \citet{shepherd2004} and \citet{yang2003} proposed to search all possible combinations of toll links using the genetic algorithm. \citet{ekstrom2012} formally represented the choice of toll links using binary variables, thus creating a mixed integer nonlinear bilevel formulation for the CBCP problem. Rather than tackling this extremely challenging formulation directly, they proposed and solved a mixed-integer linear program that functions as an approximation and provides a lower bound to the original problem. 

In summary, our reading of the literature indicates no existing algorithm promises to solve the CBCP problem at scale with a satisfactory guarantee of convergence. We intend to fill this gap by taking an approach that significantly departs from the literature.

\section{Problem setting}\label{setting}

Consider a general network $G(\cN,\cA)$ with $\cN$ being the set of nodes and $\cA$ being the set of links. We let $\cW$ denote the set of all origin-destination (OD) pairs and  $d=(d_w)_{w\in \cW}$   the vector of OD demands. The set of all the paths that travelers can choose between an OD pair  $w\in \cW$ is denoted as $\cR_w$.  Consequently, $\cR =\cup_{w\in \cW} \cR_w$ represents the set of all paths in the network. Let $v=(v_a)_{a\in \cA}$ be the vector of link flows with $v_a$ denoting the total link flow on link $a\in \cA$, and $h=(h_{r,w})_{r\in \cR_w,w\in \cW}$ be the vector of path flows with $h_{r,w}$ denoting the path flow on path $r\in \cR_w$. Feasibility dictates that link and path flows are related to each other by 
\[
v_a = \sum_{w\in \cW}\sum_{r\in \cR_w}h_{r,w}\delta_{a,r},\ \forall a\in \cA,
\]
where $\delta_{a,r}=1$ if path $r$ passes link $a$; and $0$ otherwise. Letting $\Delta = (\delta_{a,r})_{a\in \cA,r\in \cR}$ denote the link-path incidence matrix, we have $v = \Delta h$. Similarly, let $\Lambda = (\lambda_{w,r})_{w\in \cW,r\in \cR}$ denote the OD-path incidence matrix, then path flows and OD demands are related to each other by $d = \Lambda h$. The latter is also known as the flow conservation condition.

A traveler incurs a travel time, $t_a$, and a toll, $u_a$, when traversing a link $a$. The vectors for link travel times and link tolls are denoted as $t=(t_a)_{a\in \cA}$ and $u=(u_a)_{a\in \cA}$, respectively. In our setting, $u$ is independent of travelers' choice whereas $t$ is modeled as a function of link flows $v$.  Let $c =(c_{r,w})_{r\in \cR_w,w\in \cW}$ be the vector of path costs, and note that $c = \Delta^T (t + u)$.

We impose a box constraint over $u$ so that the feasible region of $u$ can be written as $U=[0,\hat{u}]$, where $\hat{u}$ is a vector of maximum permissible link tolls.  The property of $t$ is regulated by the following assumption, which is common in transportation applications.

\begin{assu}\label{ass}
The link travel time $t(v)$ is a separable, continuously differentiable, strictly increasing, and convex function of link flows $v$.
\end{assu}
For an example of $t$ that satisfies Assumption \ref{ass}, consider 
the popular Bureau of Public Roads (BPR) function, which reads,
\begin{equation}\label{bpr}
t_a(v_a) = t_{a,0}\left(1+0.15\left(\frac{v_a}{C_a}\right)^4\right),\quad \forall a\in \cA,
\end{equation}
where $t_{a,0}$ is the free-flow travel time and $C_a$ is the capacity of link $a\in \cA$. 

The travelers play a congestion game against each other, by choosing the path that gives the minimum cost. The outcome of the game is characterized by the user equilibrium (UE) conditions \citep{wardrop1952road} that state
\begin{equation} \label{eq:wadrop}
    c_{r,w} > \min_{r'\in \cR_w}\{c_{r',w}\} \Longrightarrow h_{r,w} = 0, \quad \forall w\in \cW, \  r\in \cR_w.
\end{equation}
An agent hopes to influence the travelers' decisions, hence the outcome of the game, by using its toll-setting power as a lever.  This leads to a BCP problem, whereby the agent aims to minimize the total travel time experienced by all travelers by setting the toll levels on a pre-specified set of toll links, denoted as $\bar{\cA}\subseteq \cA$ with $|\bar{\cA}|=\kappa$.  
Mathematically, the BCP problem takes the following form:
\begin{equation} \label{sbcp}
\begin{aligned}
 \rmin\limits_{{u},v}\quad & F(v)=t(v)^\top v \\ 
  \rst \quad& u\in U,\ u_a =0\ a\in \cA\backslash \bar{\cA},\\
 \quad& v\in \cS(u),
\end{aligned}
\end{equation}
where $F(v)$ represents the total travel time and $\cS(u)$ is the solution set of a traffic assignment problem
\begin{equation}\label{lp1}
\begin{aligned}
\rmin\limits_{v} \quad & f(u,v)=\sum\limits_{a\in \cA} \int_0^{v_a} (t_a(x)+u_a)dx
\\
\rst \quad & v\in \Omega\equiv\{v': v'=\triangle h,\ \Lambda h=d,\ h\ge0\},
\end{aligned}
\end{equation}
where $\Omega$ denotes the feasible region of the lower-level problem, which is defined by the flow conservation constraint, path-link relationship, and non-negativity.  The toll vector $u$ and the link flow vector $v$ represent, respectively, the upper-level and lower-level variables. Problem \eqref{lp1} is the lower-level problem parametrized by the upper-level variables $u$ and its KKT conditions are equivalent to the UE conditions \eqref{eq:wadrop} \citep{sheffi1985urban}. Thus, for each $u$ we refer to $v\in \cS(u)$ as a $u$-tolled UE link flow. 

If, in addition to setting $u$ over all $a\in \bar{\cA}$, the agent also wishes to find an optimal $\bar{\cA}$ among all subsets $\cB\subset \cA$ such that $|\bar{\cA}|\le \kappa$ (optimal in the sense the choice helps minimize the total travel time), we have to add a combinatorial optimization layer on top of the BCP problem, leading to a combinatorial BCP, or CBCP problem. Without loss of generality, we assume here that collecting toll on any link  $a\in \cA$ is feasible. If needed, we can replace $\cA$ with any of its subset. In the CBCP problem, the first-tier decision is to choose $\bar{\cA}$, of which there are $m!/(\kappa!(m-\kappa)!)$ possibilities.  Considering that the BCP problem is already intractable, solving an overwhelming number of them to address the CBCP problem is clearly not an appealing option.  To break this impasse, the next section proposes a new formulation that utilizes a cardinality function.

\section{Cardinality-constrained formulation}\label{model}

A cardinality function maps a vector $a$ to its number of non-zero components, written as  $|{\rm supp}(a)|=|\{i: a_i\not=0\}|$. Using such a function, we can write the requirement that the toll link set $\bar{\cA}$ have no more than $\kappa$ elements simply as $|{\rm supp}(u)| \le \kappa$, where $u$ is the toll vector.  This enables us to formulate the CBCP problem still as a bilevel program:
\begin{subequations}\label{p1}
\begin{align}
\rmin\limits_u &\quad  F(v) \tag{\ref{p1}{a}}\\
\rst &\quad u\in U,\tag{\ref{p1}{b}}\label{p1-a} \\
           &\quad |{\rm supp}(u)|\le \kappa,\tag{\ref{p1}{c}} \label{p1-b}\\
           &\quad v\in \cS(u). \tag{\ref{p1}{d}}\label{p1-c}
\end{align}
\end{subequations}
In contrast to Problem \eqref{sbcp}, the equality constraints $u_a =0,\ a\in \cA\backslash \bar{\cA}$ are replaced by an inequality constraint \eqref{p1-b}. By limiting the maximum number of toll links, the new formulation avoids adding another optimization layer while simultaneously exploring all possible combinations of toll locations and the associated toll levels.   Although the cardinality constraint simplifies the structure of the CBCP problem, it is a hard constraint because $|{\rm supp}(u)|$ is not continuous at points with zero components.  This discontinuity poses a significant challenge to algorithm development, which will be resolved in Section \ref{sec:PD}.  

We next transform Problem \eqref{p1} into a single-level problem, using the concept of value function. Let $\cV(u)$ denote the marginal value function of the lower-level problem, which is defined by
\begin{equation*}
\cV(u) = \inf_{v} \{ f(u,v): v\in \Omega\}.
\end{equation*}
Then, the solution set $\cS(u)$ can be expressed as
\begin{equation}\label{reformset}
\cS(u) = \{v\in \Omega: f(u,v) - \cV(u) \le 0\}.
\end{equation}
We call $f(u,v) - \cV(u)$  the \textit{gap function} of the lower-level problem under the toll level $u$, and note that the gap function is nonnegative (per definition of the marginal value function) and attains the value of zero only when the lower-level solution $v$ is at the UE state, or $v\in\cS(u)$. 

Using both the cardinality constraint and the gap function, we can now write the CBCP problem as a single-level nonlinear program:
\begin{subequations}\label{pp1}
\begin{align}
\rmin\limits_{u,v} &\quad  F(v) \tag{\ref{pp1}{a}}\\
\rst &\quad u\in U \cap U_{\kappa},\ v\in  \Omega,\tag{\ref{pp1}{b}}\label{pp1-b}\\
           &\quad f(u,v) - \cV(u)\le0,\tag{\ref{pp1}{c}}\label{pp1-c}
\end{align}
\end{subequations}
where $U_{\kappa} = \{u:|{\rm supp}(u)|\le \kappa\}$, and Constraint \eqref{pp1-c} ensures the UE conditions are always satisfied. Single-level as it may be, solving Problem \eqref{pp1} remains a difficult task because the formulation contains two intractable functions: the implicit value function $\cV(u)$ and the non-convex, discontinuous cardinality function $|{\rm supp}(u)|$. %

Before setting out to develop an efficient algorithm for solving Problem \eqref{pp1}, we first recall a few properties of the value function $\cV(u)$ that will facilitate algorithm development. 

\begin{prop}\label{prop}
Given Assumption \ref{ass}, we have the following properties: (i) The lower-level solution set $\cS(u)$ is a singleton for any $u$, and the value function $\cV(u)$ is continuously differentiable with
\begin{equation}\label{diff1}
\nabla \cV(u) = \cS(u).
\end{equation}

(ii) The value function $\cV(u)$ is concave and the gap function $f(u,v) - \cV(u)$ is block-wise convex, i.e., it is convex with respect to the upper-level variables $u$  given $v$ and convex with respect to the lower-level variables $v$ given $u$.
\end{prop}

\proof{Proof.} See Appendix \ref{appen}. \hfill

\section{Algorithm}\label{sec:PD}

In this section, we develop an efficient algorithm for solving the CBCP problem with the guarantee of convergence to a stationary point. The presentation is structured as follows. In Section \ref{sec:alg-ref}, the CBCP problem \eqref{pp1} is further transformed into a decomposable form. Section \ref{sec:alg-des} then gives an algorithm that consists of an outer and an inner loop: the latter solves a penalized approximation of the reformulation and the former adjusts the penalty factor.   Finally, Section \ref{sec:alg-conv} provides the convergence result.

\subsection{Reformulation} \label{sec:alg-ref}
As noted earlier, the cardinality constraint \eqref{p1-b} is intractable due to discontinuity and non-convexity. We seek to bypass this intractability by first replacing the toll vector $u$ with an auxiliary vector not constrained by cardinality, calling it $z$, and then forcing $z = u$.  Accordingly, we rewrite Problem \eqref{pp1} as follows:
\begin{subequations}\label{ppp1}
\begin{align}
\rmin\limits_{u,z,v} &\quad  F(v) \tag{\ref{ppp1}{a}}\\
\rst &\quad (u,v)\in U_\kappa\times\Omega,\ z\in U,\tag{\ref{ppp1}{b}}\label{ppp1-b}\\
           &\quad f(z,v) - \cV(z)\le0,\ u-z=0\tag{\ref{ppp1}{c}}\label{ppp1-c}.
\end{align}
\end{subequations}
Problem \eqref{ppp1} is more desirable because it has two block-separable constraints \eqref{ppp1-b} with two variable blocks $(u,v)$ and $z$,  each corresponding to its own feasible set, $U_\kappa \times \Omega$ and $U$, respectively.
 Problem \eqref{ppp1} still has two coupling (hence harder) constraints: $f(z,v) - \cV(z)\le0$ and $u-z=0$.  Yet, we can eliminate them by adding penalty terms into the objective function, hence obtaining a penalty approximation (PA) of Problem \eqref{ppp1} that reads
 \begin{equation}\label{pp1p}
{(\rm PA_\rho)}~~~~~
\begin{aligned}
		\rmin\limits_{u,z,v} \quad &  \Phi_\rho(u,z,v)=F(v) + \rho_1 \left(f(z,v)-\cV(z)\right)+ \rho_2\|u-z\|^2\\
		{\rst} \quad & (u,v)\in U_\kappa \times \Omega,\ z\in U,
\end{aligned}
\end{equation}
where $\rho=(\rho_1,\rho_2)>0$ is a penalty vector. 

For a given $\rho$, Problem ${(\rm PA_\rho)}$ is tractable, since it can be decomposed into two subproblems to be solved iteratively: minimizing $\Phi_\rho$ by choosing $z$ given $u$ and $v$, and minimizing $\Phi_\rho$ by setting $u$ and $v$ given $z$. For the former, we note that ${(\rm PA_\rho)}$ is convex with respect to $z$ as per Proposition \ref{prop}(ii).  The latter problem can be further decomposed since $u$ and $v$ are not interacting with each other in Problem ${(\rm PA_\rho)}$ once $z$ is fixed.  Furthermore, under Assumption \ref{ass}, Problem ${(\rm PA_\rho)}$ is convex with respect to $v$, which makes it easy to find $v$ that minimizes $\Phi_\rho$, given $z$ and independent of $u$. This leaves us with the last problem of minimizing $\Phi_\rho$ by selecting $u$ given $z$. At first glance, this seems an intractable problem in its own right, as $U_\kappa$ is non-convex. However, the part of $\Phi_\rho$ relevant to the selection of $u$ is the scaled Euclidean distance between $u$ and $z$,  i.e., $\rho_2\|u-z\|^2$.  Minimizing this distance amounts to projecting $z$ onto $U_\kappa$. As shown in the next result, this projection can be performed very efficiently, despite the seemingly intractable geometry of $U_\kappa$.

\begin{prop}\label{global}
Consider the following optimization problem:
\begin{equation}\label{projp}
\rmin_u\ \|z-u\|^2\quad {\rst}\ u\in U_\kappa.
\end{equation}
Let $I$ be the index set corresponding to the $\kappa$ largest values of $\{|z_a|:a\in \cA\}$ and $I_c$ be the complement of $I$ with respect to $\cA$. The global solution of Problem \eqref{projp}, denoted as $u^*=(u_a^*)_{a\in \cA}$, is given by 
\begin{equation} \label{eq:tollprojection}
u_a^* = \left\{\begin{array}{ll} z_a\ \ & {\rm if}\ a\in I,\\[5pt]
0\ \ & {\rm if}\ a\in I_c. \end{array}\right.
\end{equation}
\end{prop}
\proof{Proof.} See Appendix \ref{appen}.

\subsection{Description} \label{sec:alg-des}
We are now ready to present the algorithm, which consists of an inner loop and an outer loop. To the inner loop, we propose the following block coordinate descent (BCD) algorithm that solves the penalty approximation problem ${(\rm PA_\rho)}$. 

{
\setlist[itemize]{leftmargin=4.5em}
\begin{alg}\label{alg1}
Given $\rho$, find $u^*_\rho$, $v^*_\rho$, and $z^*_\rho$ that solves Problem ${(\rm PA_\rho)}$.

\begin{itemize}
\item[Step (0):] Choose an initial toll vector $z^1\in U$. Set the iteration index $r=0$.
\item[Step (1):]  Given $z^{r+1}$, solve the following minimization problem to get $(u^{r+1},v^{r+1})$:
\begin{equation}\label{alg1-sub}
\rmin\ \Phi_\rho(u,z^{r+1},v) \quad {\rst}\ u\in U_\kappa, v\in\Omega.
\end{equation}

\item[Step (2):] If $(u^{r+1},z^{r+1},v^{r+1})$ is an approximate KKT stationary solution, then stop and set $(u^*_\rho, z^*_\rho, v^*_\rho) = (u^{r+1},z^{r+1},v^{r+1})$. Otherwise, set $r=r+1$ and go to Step (3).

\item[Step (3):] For fixed $(u^r,v^r)$, solve the following convex minimization problem to get $z^{r+1}$ and go to Step (1):
    \begin{equation}\label{alg1-sub3}
      \rmin\ \Phi_\rho(u^r,z,v^r)\quad {\rst}\ z\in U.
    \end{equation}

\end{itemize}
\end{alg}
}

We proceed to discuss how Problems \eqref{alg1-sub} and \eqref{alg1-sub3} may be solved. 
Since it is partly constrained by $U_\kappa$, Problem  \eqref{alg1-sub} is non-convex. However, because $u$ and $v$ are separable in both the objective function and the constraints, Problem  \eqref{alg1-sub} can be solved by separately solving the following two problems:
\begin{equation}\label{d1}
\rmin_u\ \|u-z^{r+1}\|^2\ \ {\rst} \ \ u\in U_\kappa,
\end{equation}
and
\begin{equation}\label{d2}
\rmin_v\  F(v) + \rho_1 f(z^{r+1},v)\ \ {\rst} \ \ v\in \Omega.
\end{equation}
Per Proposition \ref{global}, Problem \eqref{d1} has a closed-form solution as given by Equation \eqref{eq:tollprojection}.
Problem \eqref{d2}, on the other hand, is a convex program with the weighted sum of the upper- and lower-level objectives as its objective, subject to the constraints of the lower-level problem. A moment of reflection suggests that it can be converted to a standard traffic assignment problem by interpreting the sum of the two terms as the potential function corresponding to a new link travel cost function 
\[
\hat{t}_a(v_a) = (1+\rho_1) t_a(v_a) + \rho_1z_a^{r+1} + v_a dt_a/dv_a.
\] 
As a result, many efficient algorithms exist that can solve Problem \eqref{d2} efficiently even over very large networks \cite[see e.g.,][for recent examples]{dial2006path,bar2010traffic,xie2016new,xie2018greedy}. Appendix \ref{appen-b} describes the algorithm used in our numerical experiments, which is based on  \cite{xie2018greedy}.

Problem \eqref{alg1-sub3} can be specified as follows:
\begin{equation}\label{form:updatez}
\rmin_z \ \rho_1(f(z,v^r)-\cV(z))+ \rho_2 \|z-u^r\|^2 \ \ {\rst} \ \ z\in U.
\end{equation}
Invoking Proposition \ref{prop} and noting the strong convexity of the quadratic term, we can easily show that Problem \eqref{form:updatez} is a strongly convex program with a simple box constraint (which defines $U$). Hence, it can be solved to global optimality using gradient descent methods.  We devise and present such an algorithm in Appendix \ref{appen-c}.

A solution obtained by Algorithm \ref{alg1} for a given penalty vector $\rho$ may not be accepted as the solution to the original Problem \eqref{ppp1}, since Constraint \eqref{ppp1-c} is not strictly enforced in the penalty approximation problem ${(\rm PA_\rho)}$.  If the solution violates these constraints, we need to adjust the penalty vector and re-solve   ${(\rm PA_\rho)}$, until the violations are eliminated. This forms the outer-loop component of the proposed algorithm, referred to as the {{penalized block coordinate descent (PBCD) algorithm}}.

{
\setlist[itemize]{leftmargin=4.5em}
\begin{alg}\label{mainalg1}
Given $G(\cN,\cA)$, link performance function $t(v)$, and the maximum number of toll links $\kappa$, find $u^*$, $v^*$ and $z^*$ that solve Problem \eqref{ppp1}.
\begin{itemize}
\item[Step (0):] Choose $z_0^1\in U$, $\rho^{1}=(\rho_1^{1},\rho_2^{1})>0$, $\gamma=(\gamma_1,\gamma_2) > 1$, a convergence criterion $\varepsilon>0$, a feasible solution $(u^{feas},v^{feas})$ to Problem \eqref{p1}, and a positive constant $\Upsilon$ such that
\[
\Upsilon \ge \max\{F(v^{feas}),\min\limits_{u\in U_\kappa,v\in\Omega}\Phi_{\rho^{1}}(u,z^1_0,v)\}.
\]
Set $k=1$.
\item[Step (1):] Solve Problem ${(\rm PA_{\rho^{k}})}$ by Algorithm \ref{alg1} to get $(u^k,z^k,v^k)$, using  $z_0^k$ as the initial toll vector.
\item[Step (2):] If $f(z^k,v^k)-\cV(z^k) \le \varepsilon$ and $\|u^k-z^k\|\le \varepsilon$,  stop; otherwise, set $\rho^{k+1} = (\gamma_1\rho_1^{k},\gamma_2\rho_2^{k})$.

\item[Step (3):] Set $z_0^{k+1}=z^k$ if 
$$\min\limits_{u\in U_\kappa,v\in\Omega}\Phi_{\rho^{k+1}}(u,z^k,v)\le \Upsilon;$$ otherwise, set $z_0^{k+1} = u^{feas}$. Set $k=k+1$ and return to Step (1).

\end{itemize}
\end{alg}
}

The  feasible solution $(u^{feas},v^{feas})$ to Problem \eqref{p1} can be obtained by choosing an arbitrary $u^{feas}\in U\cap U_\kappa$ and then setting $v^{feas}$ as  the $u^{feas}$-tolled user equilibrium solution. When solving the penalty approximation problem, the constant $\Upsilon$ guides the choice of the initial toll vector such that the generated sequence $\{\Phi_{\rho^{k}}\}_{k=1}^\infty$ is bounded above by a constant regardless of the magnitude of the penalty factors.  As we shall see later, this feature ensures the feasibility can be achieved when the penalty factor approaches infinity, hence the key to the convergence of Algorithm \ref{mainalg1}. %

\subsection{Convergence} 
\label{sec:alg-conv}

We establish the convergence results for Algorithms \ref{alg1} and \ref{mainalg1} separably in this section. Since Problem ${(\rm PA_\rho)}$ is non-convex, the best that Algorithm \ref{alg1} can be expected to achieve is a stationary, or KKT stationary point  \citep[e.g.,][]{treiman1999,nocedal1999}. We begin by formally defining KKT points.

Following \citet[Chapter 6]{VA}, we define the KKT point of an optimization problem as a point where the negative gradient of its objective function lies in the limiting normal cone to the feasible set at that point. Let $W$ be a closed set and $\bw\in W$. The limiting normal cone to $W$ at $\bw$ is defined as
\[
N_W(\bw) =\{z: \exists w^k\to \bw, z^k \to z \ {\rm with}\ (z^k)^T(w-w^k)\le o(\|w-w^k\|)\ \forall w\in W\}.
\]
When $W$ is convex, the limiting normal cone $N_W(\bw)$ coincides with the classical normal cone in convex analysis, defined as 
\[
N^c_W(\bw)=\{z:z^T(w-\bw)\le 0\ \forall w\in W\}.
\]

\begin{defi}\label{defi}
Consider a minimization problem with two block variables $x$ and $y$ as follows:
\begin{equation}\label{eq:canonical}
\begin{aligned} 
\rmin &\quad f(x,y) \nonumber\\
\rst &\quad g(x,y)\le0,\ h(x,y)=0,\\
           &\quad x\in X,\ y\in Y.\nonumber
\end{aligned}
\end{equation}
We say that a feasible point $(\bx,\by)$ is a KKT stationary point if there exist multipliers $\lambda$ and $\mu$ such that
\begin{equation}
\begin{aligned}
& 0\in \nabla f(\bx,\by) + \nabla g(\bx,\by) \lambda + \nabla h(\bx,\by)\mu + N_{X}(\bx)\times N_{Y}(\by),\\
& \lambda \ge0,\ g(\bx,\by)^\top \lambda=0.
\end{aligned}
\end{equation}
\end{defi}

Problem ${(\rm PA_\rho)}$ is of the canonical form given by \eqref{eq:canonical}.  The following result asserts that Algorithm \ref{alg1} is a strictly descent algorithm and converges to a KKT stationary solution of Problem ${(\rm PA_\rho)}$.

\begin{thm}\label{sub-con1}
Assume that $\{(u^r,z^r,v^r)\}_{r=1}^\infty$ is an infinite sequence generated by Algorithm \ref{alg1} where $u^r$ is derived as in Proposition \ref{global}.
\begin{itemize}

\item[(i)] Algorithm \ref{alg1} is a strictly descent algorithm, that is,
\begin{equation}\label{monoto}
\Phi_\rho(u^{r+1},z^{r+1},v^{r+1}) < \Phi_\rho(u^r,z^r,v^r) \le \min\limits_{u\in U_\kappa,v\in\Omega}\Phi_{\rho}(u,z^1,v),\quad \forall r\ge1.
\end{equation}

\item[(ii)] The sequence $\{(u^r,z^r,v^r)\}_{r=1}^\infty$ is bounded and $\{\Phi_\rho(u^r,z^r,v^r)\}_{r=1}^\infty$ has a unique limit. Any accumulation point $(u^*,z^*,v^*)$ of $\{(u^r,z^r,v^r)\}_{r=1}^\infty$ is a KKT stationary solution of ${(\rm PA_{\rho})}$.
\end{itemize}
\end{thm}
\proof{Proof.} See Appendix \ref{appen}.

To establish the convergence of Algorithm \ref{mainalg1} (PBCD), we focus on the number of iterations required to obtain a sufficiently precise approximate solution to Problem  \eqref{pp1}. We begin with a special case where the penalty approximation problem can be solved globally and hence a stronger convergence result is secured.  
\begin{thm}\label{mainthm1}
Let $\{(u^k,z^k,v^k)\}_{k=1}^\infty$ be a sequence generated by Algorithm \ref{mainalg1} and suppose $(u^k,z^k,v^k)$ is a globally optimal solution of ${(\rm PA_{\rho_k})}$ for each $k$. When 
\begin{equation}\label{number1}
k\ge\max\left(\frac{\ln(F^*-F^l)-\ln(\varepsilon \rho^1_1)}{\ln \gamma_1}+1,\frac{\ln(F^*-F^l)-\ln(\varepsilon^2 \rho^1_2)}{\ln \gamma_2}+1\right),
\end{equation}
where $\gamma_i, \rho_i^1, i = 1,2,$ are parameters used by Algorithm \ref{mainalg1}, $F^*=\min\limits_{v\in \cS(u),u\in U\cap U_\kappa} F(v)$ and $\ F^l = \min\limits_{v\in \Omega} F(v)$, the derived solution $(u^k,z^k,v^k)\in U_\kappa\times U\times \Omega$ is an $\varepsilon$-optimal solution to Problem \eqref{pp1}, in the sense that it satisfies
\begin{equation}\label{appopt}
F(v^k)\le F^*,\ f(z^k,v^k) - \cV(z^k)\le \varepsilon, \ \|u^k-z^k\|\le \varepsilon.
\end{equation}
\end{thm}
\proof{Proof.} See Appendix \ref{appen}.

When $\varepsilon$ in Condition \eqref{appopt} reaches zero for a given $k$,  $u^k \in U_\kappa$ implies that we have identified a tolling strategy with at most $\kappa$ toll links. The last two inequalities with $\varepsilon=0$ in Condition \eqref{appopt} imply that $v^k$ is a $u^k$-tolled user equilibrium link flow. This and the first inequality imply that $(u^k,v^k)$ is a globally optimal solution to Problem \eqref{pp1}. For a positive $\varepsilon$, the constraints are slightly relaxed and thus a better system objective may be achieved, as indicated in Condition \eqref{appopt}.

If the penalty approximation problem ${(\rm PA_{\rho})}$ cannot be solved globally  --- which is more line with the problem considered herein  ---  the convergence result is weakened to the following.
\begin{thm}\label{mainthm}
Let $\{(u^k,z^k,v^k)\}_{k=1}^\infty$ be a sequence generated by Algorithm \ref{mainalg1} (PBCD). When the iteration number satisfies
\begin{equation}\label{number}
k\ge\max\left(\frac{\ln(\Upsilon-F^l)-\ln(\varepsilon \rho_1^1)}{\ln \gamma_1}+1,\frac{\ln(\Upsilon-F^l)-\ln(\varepsilon^2 \rho_2^1)}{\ln \gamma_2}+1\right),
\end{equation}
where $\Upsilon,\gamma_i,\rho_i^1, i = 1,2,$ are parameters given in Algorithm \ref{mainalg1} and $F^l$ is given in Theorem \ref{mainthm1}.
The derived solution $(u^k,z^k,v^k)\in U_\kappa\times U\times \Omega$ is an $\varepsilon$-approximate KKT stationary solution of Problem \eqref{ppp1}, in the sense that it satisfies
\begin{equation}\label{thm-res}
\begin{aligned}
& 0\in \nabla_z f(z^k,v^k) -\nabla \cV(z^k) + N_{U_\kappa}(u^k) + N_U(z^k),\\
& 0\in \nabla F(v^k) + \rho_k \nabla_v f(z^k,v^k) + N_\Omega(v^k),\\
& f(z^k,v^k)-\cV(z^k) \le \varepsilon,\ \|u^k-z^k\| \le \varepsilon.
\end{aligned}
\end{equation}
\end{thm}
\proof{Proof.} See Appendix \ref{appen}.

In other words, instead of getting a sufficiently precise global solution, the proposed PBCD algorithm only guarantees a solution sufficiently close to a KKT point given enough iterations. 

\section{Numerical study}\label{sec:num}

\begin{figure}[ht]
    \centering
    \includegraphics[width=0.6\textwidth]{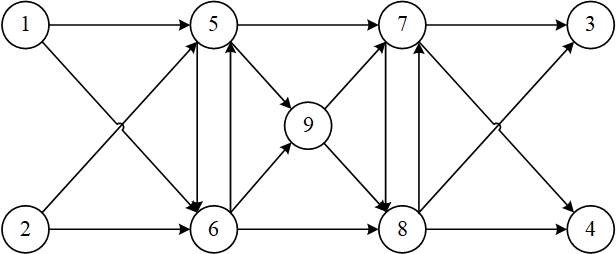} 
    \centering
    \caption{Topology of the Hearn's network.}
    \label{fig:nn}
\end{figure}

To examine the performance of the proposed PBCD algorithm (i.e. Algorithm \ref{mainalg1}), we conduct in this section a set of numerical experiments on three networks frequently used in the transportation literature as benchmarks: the network from \cite{hearn1998solving}, hereafter referred to as the Hearn's network  (Section \ref{sec:hearn}), the Sioux-Falls network (Section \ref{sec:realistic-1}), and the Chicago-Sketch network (Section \ref{sec:realistic-2}). 
The topology of the Hearn's network is shown in Figure \ref{fig:nn}. For the other two, Sioux-Falls has 24 nodes, 76 links, and 528 OD pairs, and  Chicago-Sketch has 933 nodes, 2,950 links, and 93,513 OD pairs.
In all networks, the travel time function takes the BPR form \eqref{bpr}.
For more details of the two larger networks, the reader may consult the Transportation Networks GitHub Repository \citep{networks-github}.

On each network, we run PBCD in various settings.
In each test,  the effectiveness of a tolling scheme obtained by our algorithm, denoted as $u^*$, is evaluated as follows. We first compute the UE link flow pattern $v^*$ under the tolling scheme $u^*$. Then, $F^* = F(v^*)$, i.e., the total travel time induced by $u^*$, is compared with two reference points: (1) $F^{so} = F(v^{so})$, i.e., the total travel time at the system optimal (SO) link flow pattern $v^{so}$; (2) $F^{ue} = F(v^{ue})$, i.e., the total travel time at the no-toll UE link flow pattern $v^{ue}$.  The effectiveness of $u^*$ is gauged by the ``relative excessive delay" (R.E.D.) at $v^*$, computed by
\begin{equation}
    \frac{F^* - F^{so}}{F^{ue} - F^{so}}.
\end{equation}
Clearly, the relative excessive delay must range between 0 and 1, and the closer to 0, the better.

{The performance of the algorithm is affected by several key hyperparameters, including the initial toll levels $z_0^1$, the termination criteria in Step (2), and the penalty adjustment factors $\gamma_1$ and $\gamma_2$ in Algorithm \ref{mainalg1}.  However, our preliminary experiments indicate that achieving a satisfactory performance consistently across different instances of CBCP problems does not require retuning most parameters. In our experiments, we always set $\gamma_1=1.8$ and $\gamma_2 =5.0$.  As for the convergence test in Algorithm \ref{mainalg1}, we adopt the following termination conditions in Step 2:
\begin{equation*}
    \frac{f(z^k,v^k)-\cV(z^k)}{\max\{f(z^k,v^k),1\}}\le\varepsilon_1,\quad
    \frac{\|u^k-z^k\|}{\max\{\|u^k\|,1\}}\le\varepsilon_2,
\end{equation*}
where $\varepsilon_1 = 0.0001, \varepsilon_2 = 0.001$ in all experiments. 
We opted for the use of relative errors rather than absolute errors because the latter provides smoother convergence by dynamically adapting to the proper scale.  In our experience, the only hyperparameter that needs retuning is the initial toll value $z^1_0$. For this parameter, we suggest the following rule of thumb which was adopted in our experiments based on trial-and-error: $z^1_0$ is set to zero for large $\kappa$ (e.g. when $\kappa> 0.2 \cdot |\cA|$) and to a nonzero value (e.g., $z^1_0 = \mathbf{1.0}$) otherwise.}

Finally, whenever a UE problem needs to be solved, the improved greedy algorithm by \citet{xie2018greedy}, an efficient path-based UE algorithm, is employed.

The algorithm is coded using the toolkit of network modeling, a C++ class library specialized in modeling transportation networks \citep{nie2006programmer}. Unless otherwise specified, all numerical results reported in this section were produced based on C++ on a Windows 11 64-bit laptop with Intel(R) 11th Gen CPU i7-11800H 2.30GHz and 16G RAM.

\subsection{Hearn's network}
\label{sec:hearn}

Table \ref{tab:tab1} reports the capacity $C_a$ and free-flow travel time $t_{a, 0}$ of all links in Hearn's network, as well as two solutions as reference: the SO solution \(v^{so}\) and the no-toll UE solution \(v^{ue}\). The total travel times at the SO and UE solutions are \(F^{so} = 37.57\) and \(F^{ue} = 40.93\), respectively, which are given in the last row of Table \ref{tab:tab1}. The difference of 3.36 represents the maximum potential reduction in total travel time that can be achieved through the application of congestion pricing.

\begin{table}[ht]
\centering
\footnotesize
\caption{Network data, SO link flows, UE link flows, and MTL toll levels for Hearn's network.}
\label{tab:tab1}
\setlength{\tabcolsep}{10 mm}{
\begin{tabular}{ccccc} \toprule
Link & $C_a$ & $t_{a, 0}$ & $v^{so}$ & $v^{ue}$  \\
\midrule
1-5 & 12 & 5 &  9.41 & 8.16 \\
1-6 & 18 & 6 & 20.59 & 21.84  \\
2-5 & 35 & 3 & 38.33 & 47.37 \\
2-6 & 35 & 9 &  31.67& 22.63 \\
5-6 & 20& 9 &  0.00 & 0.00 \\
5-7 & 11& 2 &  21.30 & 27.84  \\
5-9 & 26& 8 & 26.44 & 27.69 \\
6-5 & 11 & 4 &0.00 & 0.00 \\
6-8 & 33 & 6 &  39.47 & 44.47 \\
6-9 & 32 & 7 &  12.78 & 0.00 \\
7-3 & 25 & 3 &  29.61 & 38.16 \\
7-4 & 24 & 6 & 20.76 & 17.37 \\
7-8 & 19 & 2 &  0.00 & 0.00 \\
8-3 & 39 & 8 &  10.39 & 1.84 \\
8-4 & 43 & 6 & 39.24 & 42.63  \\
8-7 & 36 & 4 &  0.00 & 0.00 \\
9-7 & 26 & 4 & 29.06 & 27.69 \\
9-8 & 30 & 8 & 10.16 & 0.00 \\ \midrule
Total travel time & - & - & 37.57& 40.93\\
\bottomrule
\end{tabular}
}
\end{table}

In Hearn's network, it is not necessary to impose tolls on all links to bring UE to SO. Instead, according to \citet{hearn1998solving},  charging tolls  on {as few as} five links suffices. This tolling scheme, referred to as the minimum-toll-location (MTL) solution hereafter, is shown in the last column of Table \ref{tab:hearn-solution}. For this reason, we will only test the cases when $\kappa = 1, \ldots, 5$. For each $\kappa$, we will first run PBCD for solving Problem \eqref{p1} and then compare the derived solution with a globally optimal solution to Problem \eqref{p1}, which is obtained by a brute-force approach, briefly described below. {Given the size of the network, enumerating {all combinations for toll links} is practical (for $\kappa \leq 5$, there are mere 8,568 combinations). For each combination, toll links are set, hence Problem \eqref{p1} is reduced to a standard bilevel congestion pricing (BCP) problem. For each BCP problem, we search for a global solution as follows.
\begin{itemize}
    \item A ``coarse" grid search is performed first, in which the objective function value is evaluated at each point of a multidimensional grid of link toll levels.
    \item The best solution from the grid search is then used as the initial solution to Powell's conjugate direction method (referred to as Powell's method hereafter), a convenient local search algorithm that does not rely on first-order information. To provide an objective benchmark, we employ an implementation of Powell's method in Python's \texttt{scipy.optimize} package.
    \item The solution generated by Powell's method is taken as a sufficiently good global solution to the BCP problem.
\end{itemize} 
After all BCP problems are solved using the above procedure,  the best solution, along with the corresponding toll links, is accepted as an optimal solution to Problem \eqref{p1}.
}

The solutions generated by PBCD and the brute-force global search with $\kappa = 1, \ldots, 5$ are reported side by side
in Table \ref{tab:hearn-solution}.  PBCD consistently achieves the global solution in all cases. When $\kappa = 5$, the solution also coincides with the MTL solution, which offers another assurance that the solution obtained is indeed globally optimal. In one case ($\kappa  =4$), PBCD apparently identifies an alternative global solution: the same total travel time as obtained by the global search but a rather different tolling scheme.  This result positively confirms that PBCD can simultaneously identify the optimal set of toll links and determine the optimal toll levels.

\setlength{\tabcolsep}{3.5pt}
\begin{table}[ht]
  \centering
  \caption{Tolling schemes given by PBCD (Algorithm \ref{mainalg1}), a global search, and the minimum-toll-location (MTL) solution from \citet{hearn1998solving}.}
  \footnotesize
    \begin{tabular}{ccccccccccccccccc}
    \toprule
    \multirow{2}[4]{*}{Link } & \multicolumn{2}{c}{$\kappa = 1$} &       & \multicolumn{2}{c}{$\kappa = 2$} &       & \multicolumn{2}{c}{$\kappa = 3$} &       & \multicolumn{2}{c}{$\kappa = 4$} &       & \multicolumn{2}{c}{$\kappa = 5$} &       & \multirow{2}[4]{*}{MTL} \\
\cmidrule{2-3}\cmidrule{5-6}\cmidrule{8-9}\cmidrule{11-12}\cmidrule{14-15}          & PBCD & Global &       & PBCD & Global &       & PBCD & Global &       & PBCD & Global &       & PBCD & Global &       &  \\
    \midrule
    1-5   & 0     & 0     &       & 0     & 0     &       & 0     & 0     &       & 0     & 0     &       & 0     & 0     &       & 0 \\
    1-6   & 0     & 0     &       & 0     & 0     &       & 0     & 0     &       & 0     & 0     &       & 0     & 0     &       & 0 \\
    2-5   & 0     & 0     &       & 0     & 0     &       & \textcolor[rgb]{ 1,  0,  0}{4.00} & \textcolor[rgb]{ 1,  0,  0}{4.00} &       & \textcolor[rgb]{ 1,  0,  0}{4.00} & \textcolor[rgb]{ 1,  0,  0}{4.00} &       & \textcolor[rgb]{ 1,  0,  0}{4.00} & \textcolor[rgb]{ 1,  0,  0}{4.00} &       & \textcolor[rgb]{ 1,  0,  0}{4.00} \\
    2-6   & 0     & 0     &       & 0     & 0     &       & 0     & 0     &       & 0     & 0     &       & 0     & 0     &       & 0 \\
    5-6   & 0     & 0     &       & 0     & 0     &       & 0     & 0     &       & 0     & 0     &       & 0     & 0     &       & 0 \\
    5-7   & \textcolor[rgb]{ 1,  0,  0}{8.00} & \textcolor[rgb]{ 1,  0,  0}{8.00} &       & \textcolor[rgb]{ 1,  0,  0}{8.00} & \textcolor[rgb]{ 1,  0,  0}{8.00} &       & \textcolor[rgb]{ 1,  0,  0}{8.00} & \textcolor[rgb]{ 1,  0,  0}{8.00} &       & \textcolor[rgb]{ 1,  0,  0}{8.00} & \textcolor[rgb]{ 1,  0,  0}{8.00} &       & \textcolor[rgb]{ 1,  0,  0}{11.20} & \textcolor[rgb]{ 1,  0,  0}{11.20} &       & \textcolor[rgb]{ 1,  0,  0}{11.20} \\
    5-9   & 0     & 0     &       & 0     & 0     &       & 0     & 0     &       & 0     & 0     &       & 0     & 0     &       & 0 \\
    6-5   & 0     & 0     &       & 0     & 0     &       & 0     & 0     &       & 0     & 0     &       & 0     & 0     &       & 0 \\
    6-8   & 0     & 0     &       & 0     & 0     &       & 0     & 0     &       & 0     & 0     &       & \textcolor[rgb]{ 1,  0,  0}{7.20} & \textcolor[rgb]{ 1,  0,  0}{7.20} &       & \textcolor[rgb]{ 1,  0,  0}{7.20} \\
    6-9   & 0     & 0     &       & 0     & 0     &       & 0     & 0     &       & 0     & 0     &       & 0     & 0     &       & 0 \\
    7-3   & 0     & 0     &       & 0     & 0     &       & 0     & 0     &       & \textcolor[rgb]{ 1,  0,  0}{0.02} & 0     &       & \textcolor[rgb]{ 1,  0,  0}{4.00} & \textcolor[rgb]{ 1,  0,  0}{4.00} &       & \textcolor[rgb]{ 1,  0,  0}{4.00} \\
    7-4   & 0     & 0     &       & 0     & 0     &       & 0     & 0     &       & 0     & \textcolor[rgb]{ 1,  0,  0}{7.47} &       & 0     & 0     &       & 0 \\
    7-8   & 0     & 0     &       & 0     & 0     &       & 0     & 0     &       & 0     & 0     &       & 0     & 0     &       & 0 \\
    8-3   & 0     & 0     &       & 0     & 0     &       & 0     & 0     &       & 0     & 0     &       & 0     & 0     &       & 0 \\
    8-4   & 0     & 0     &       & 0     & 0     &       & \textcolor[rgb]{ 1,  0,  0}{4.00} & \textcolor[rgb]{ 1,  0,  0}{4.00} &       & \textcolor[rgb]{ 1,  0,  0}{4.00} & \textcolor[rgb]{ 1,  0,  0}{11.47} &       & 0     & 0     &       & 0 \\
    8-7   & 0     & 0     &       & 0     & 0     &       & 0     & 0     &       & 0     & 0     &       & 0     & 0     &       & 0 \\
    9-7   & 0     & 0     &       & 0     & 0     &       & 0     & 0     &       & 0     & 0     &       & \textcolor[rgb]{ 1,  0,  0}{3.20} & \textcolor[rgb]{ 1,  0,  0}{3.20} &       & \textcolor[rgb]{ 1,  0,  0}{3.20} \\
    9-8   & 0     & 0     &       & 0     & 0     &       & 0     & 0     &       & 0     & 0     &       & 0     & 0     &       & 0 \\
    \midrule
    R.E.D.   & 53.1\% & 53.1\% &       & 53.1\% & 53.1\% &       & 13.8\% & 13.8\% &       & 13.8\% & 13.8\% &       & 0.00\% & 0.00\% &       & 0.00\% \\
    \bottomrule
    \end{tabular}%
  \label{tab:hearn-solution}%
\end{table}%

The above analysis can help the toll designer make informed decisions regarding the trade-offs between alleviating congestion externality and managing the impact of tolling. A case in point is to choose a proper set of toll links.  As shown in Table \ref{tab:hearn-solution},  increasing the number of toll links from 1 to 2  and from 3 to 4 makes no difference to the total travel time at all.  Furthermore, increasing $\kappa$ from 3 to 5  brings modest improvements, but no more improvements are possible beyond that. With this knowledge, the designer only needs to compare three alternatives: $\kappa =1, 3$ and 5. %

\subsection{Sioux-Falls network}
\label{sec:realistic-1}

In the Sioux-Falls network, the total travel times at the UE and SO solutions are 124,670 and 119,904, respectively. To test the performance of our approach, PBCD is executed with $\kappa = 10, 20, \ldots, 60$. Note that the selection of $\kappa$ in this context serves only to evaluate the proximity of the algorithm's solutions to the system optimum. It does not indicate that such a high number of tolled links is desirable or realistic in practice. The results are shown in Table \ref{tab:sf}, where the second row reports the R.E.D. corresponding to the solution obtained by Algorithm \ref{mainalg1}, and the third row reports the CPU time required for the algorithm to converge. From Table \ref{tab:sf}, we observe that the computational time required for convergence is under 2 minutes for all tested $\kappa$.  A general trend is that the CPU time required  to converge decreases with the increase of $\kappa$. When $\kappa = 30$ (out of 76), the R.E.D. is around 1\%, indicating that the solution is very close to the SO solution. 
 \setlength{\tabcolsep}{4pt}
\begin{table}[ht]
\centering
\caption{Numerical results under different $\kappa$ on Sioux Falls network.}
\label{tab:sf}
\footnotesize
\begin{tabular}{cccccccccc}
\toprule
$\kappa$ & 10 & 20 & 30 & 40 & 50 & 60 \\
\midrule
R.E.D. & 25.0\%  & 6.7\% & 1.3\% & 0.02\%  & 0.00\% & 0.00\%  \\
CPU time (s) & 80.8 & 66.5 & 49.8 &  35.0 & 11.8 & 4.9 \\
\bottomrule
\end{tabular}

\end{table}

We next compare the solutions obtained by PBCD with those given by several alternative algorithms that separate the choice of toll links from the determination of toll levels. These algorithms differ from each other only on how the toll link set is set.  The goal here is to verify whether the proposed algorithm outperforms the popular heuristics in terms of identifying the most promising locations for levying tolls. 

The first alternative, referred to as PBCD' hereafter, simply takes the toll links identified by PBCD. All other algorithms utilize some heuristics (referred to as \textbf{H1}--\textbf{H4} hereafter) to determine the set of toll links.   To describe them, let us first define $\tilde{\mathcal{A}} = \{a \in \mathcal{A}: v_{a}^{ue} > v_{a}^{so}\}$, i.e., the set of links on which the UE flow is higher than the SO flow. 
\begin{itemize}
    \item \textbf{H1} and \textbf{H2} \citep{harks2015computing}  select links from $\tilde{\mathcal{A}}$  with the top $\kappa$ largest values of, respectively,
    \begin{equation*}
        t_a'(x_a^{ue}) \cdot x_a^{ue} \quad \text{and} \quad t_a'(x_a^{ue}) \cdot x_a^{ue} - t_a'(x_a^{so}) \cdot x_a^{so}. 
    \end{equation*}
    \item \textbf{H3} \citep{harks2015computing} selects the top $\kappa$ links with the largest values of $ v_a^{ue} - v_a^{so}$. 
    \item \textbf{H4} first performs a sensitivity analysis using \citet{yangandhuang2005}'s method to derive the derivative of the total travel time with respect to a link toll at a no-toll UE solution. Then the  top $\kappa$ links with the  most negative derivatives are selected.
\end{itemize}

Once the toll link set is given,  the problem is reduced to a standard BCP problem, which is subsequently solved by Powell's conjugate direction method. We start the local search with multiple initial solutions and accept the best local solution. Powell's method is implemented similarly as described in Section \ref{sec:hearn}.

\begin{figure}[ht]
\vskip 0.1in
\centering

\begin{subfigure}[b]{0.75\textwidth}
\includegraphics[width=1\columnwidth]{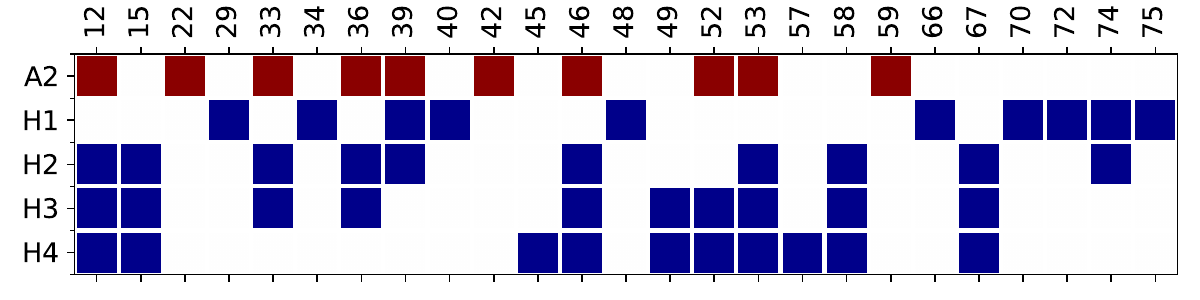}
\vspace{-1.8em}
\caption{$\kappa = 10$.}
\label{fig:sf-1}
\end{subfigure}

\vspace{1em}
\begin{subfigure}[b]{1\textwidth}
\includegraphics[width=0.98\columnwidth]{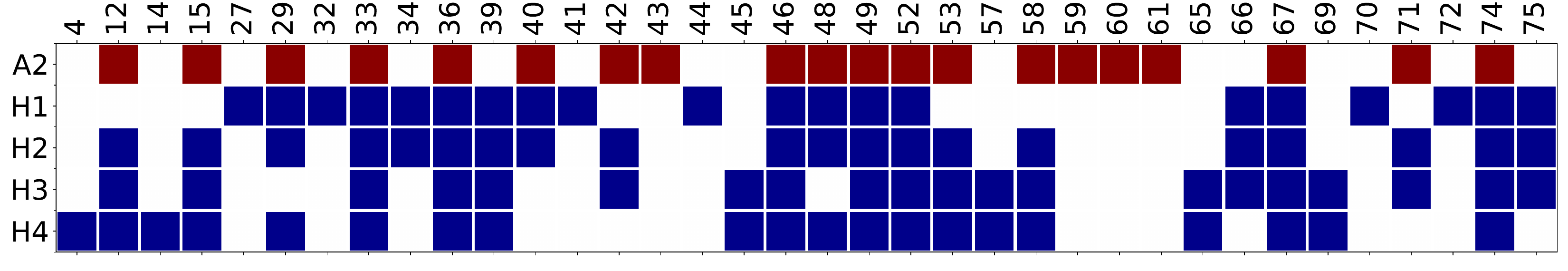}
\vspace{-0.5em}
\caption{$\kappa = 20$.}
\label{fig:sf-2}
\end{subfigure}
\caption{Toll links selected by PBCD (Algorithm \ref{mainalg1}) and \textbf{H1}--\textbf{H4}.}
\label{fig:sf}
\end{figure}

For brevity, we only make the comparison under two settings: $\kappa = 10$ and $\kappa = 20$. 
Table \ref{tab:sf-2}  shows that the solution given by PBCD is as good as that achieved by any of the five heuristics --- in most cases, it is far better. PBCD' emerges as a close second, which is expected, given that it employs what we believe to be a near-optimal toll link set.   Of the four heuristics, it appears that \textbf{H2} and \textbf{H3} outperform the other two with significant margins.  Figure \ref{fig:sf} visualizes the toll links identified by PBCD and the five heuristic schemes.  We can see that, in either scenario ($\kappa = 10$ or $\kappa = 20$), no two methods produce the same set of toll links.  This is, of course, not particularly surprising.  It is worth noting, however, that PBCD found a few links that no other method was able to uncover. Evidently, these links contribute to the better performance of PBCD.

\setlength{\tabcolsep}{2.5pt}
\begin{table}[ht]
  \centering
  \caption{Qualify of the solutions obtained by PBCD (Algorithm \ref{mainalg1}) vs those obtained using PBCD' and \textbf{H1}--\textbf{H4}.}
  \label{tab:sf-2}%
    \vspace{0.25em}
    \begin{subtable}[ht]{0.47\textwidth}
  \footnotesize
  \centering
  \caption{$\kappa = 10$.}
  \vspace{-2pt}
        \begin{tabular}{ccccccc}
    \toprule
    Algorithm & PBCD    & PBCD' &\textbf{H1}    & \textbf{H2}    & \textbf{H3}    & \textbf{H4} \\
    \midrule
    R.E.D. & 25.0\% & 25.0\% & 75.4\% & 33.0\% & 33.2\% & 45.8\% \\
    \bottomrule
    \end{tabular}%
  \label{tab:sf-2-1}%
  \end{subtable}
\hfill
    \begin{subtable}[ht]{0.47\textwidth}
    \footnotesize
    \centering
    \caption{$\kappa = 20$.}
    \vspace{-2pt}
    \begin{tabular}{ccccccc}
    \toprule
    Algorithm & PBCD   & PBCD' & \textbf{H1}    & \textbf{H2}    & \textbf{H3}    & \textbf{H4} \\
    \midrule
    R.E.D. & 6.7\% & 7.1\% & 54.2\% & 12.1\% & 17.7\% & 45.0\% \\
    \bottomrule
    \end{tabular}%
  \label{tab:sf-2-2}%
  \end{subtable}
\end{table}%

The results again highlight the effectiveness of the proposed algorithm in killing two birds with one stone.  In the literature, the  task of selecting toll links is often separated from setting toll levels because the former involves combinatorial optimization, which is considered intractable. PBCD not only integrates the two tasks seamlessly but also delivers quality solutions unmatched by conventional methods.

\subsection{Chicago-Sketch network}
\label{sec:realistic-2}

By today's standard, Chicago-Sketch is at most a medium-sized network in transportation planning practice. However, for a combinatorial optimization problem like CBCP, such a network (with nearly 3,000 links) is enormous.  In this section, we test the proposed algorithm on this network to showcase its applicability in real-world applications. We vary the value of $\kappa$ from 10 to 1600 as shown in Table \ref{tab:tab5}. The selection of $\kappa$ in this context serves only to evaluate the proximity of the algorithm's solutions to the system optimum. It does not indicate that such a high number of tolled links is desirable or realistic in practice.
\setlength{\tabcolsep}{4pt}
\begin{table}[ht]
\centering
\small
\renewcommand{\arraystretch}{1.5}
\caption{Numerical results under different $\kappa$ on Chicago-Sketch network.}
\label{tab:tab5}
\setlength{\tabcolsep}{1.5 mm}{
\begin{tabular}{ccccccccc}
\toprule
$\kappa$ & 10 & 50 & 100 & 200 & 500 & 800 & 1200 & 1600 \\
\midrule
R.E.D. & 83.7\% & 58.9\% & 44.7\% & 26.9\% & 6.8\% & 2.5\% & 1.4\% &  0.5\% \\
CPU time (s) & 1418 & 1398 & 1360 & 1235 & 1124 & 1111 & 1119 &  1265 \\ 
\bottomrule
\end{tabular}
}
\end{table}

The results show that even on such a large network, the PBCD algorithm was able to consistently solve the CBCP problem in less than 25 minutes (ranging from 18.5 to 23.7 minutes).  The value of $\kappa$ still affects the computational time, although the impact appears relatively modest. While the PBCD algorithm only promises a solution sufficiently close to a KKT point, we can see that with 1600 toll links, it reaches a solution within 0.5\% of the SO solution --- the true gap is likely much smaller than that crude estimate.  The results confirm again that adding more links to the toll set has a clear diminishing marginal return.  With 500 toll links, the PBCD algorithm has already achieved the vast majority of the potential gains (more than 93\%); and adding another 1100 offers an improvement of just 6\%. Again,  such information could play a meaningful role in aiding relevant decision-making processes.

\section{Concluding remarks}\label{sec:con}
In this study, we tackled the combinatorial bilevel congestion pricing (CBCP) problem, a variant of the mixed network design problem. This is a well-known computational challenge that, despite significant attention in the literature --- particularly in transportation --- remains unresolved to a satisfactory degree. Conventional wisdom suggests that these problems are intractable since they have to be formulated and solved with a significant number of integer variables. We showed that the CBCP problem, which aims to minimize the total system travel time by choosing both toll locations and levels, is amenable to a scalable local algorithm that guarantees convergence to an approximate KKT point.  We are able to apply the algorithm to solve, in about 20 minutes, the CBCP problem with up to 3,000 links.  To the best of our knowledge, no existing algorithm can solve CBCP problems at such a scale while providing any assurance of convergence. In small instances, our numerical experiments verified the ability of the algorithm to find the optimal toll locations and toll levels simultaneously and consistently. In larger cases for which the global optima are unknown, our algorithms were found to outperform existing heuristics with wide margins. 

Our approach is novel in that it eliminates the use of integer variables altogether, instead introducing a cardinality constraint that limits the number of toll locations to a pre-specified upper bound.  However, a bilevel program with the cardinality constraint remains a formidable challenge.  A path forward was forged by taking advantage of the fact that the projection onto the cardinality constraint is available in closed form. This enables us to transform the bilevel program into a block-separable, single-level optimization problem that can be solved efficiently after penalization and decomposition.  Importantly, we established the convergence result for the proposed PBCD algorithm, proving that, under mild conditions, it is guaranteed to reach an approximate KKT point of the original problem with sufficient precision. 

Despite the impressive performance, the PBCD algorithm is not a panacea to general mixed network design problems. It is developed using a value function-based reformulation, which relies on the separability of link travel times, as highlighted in Assumption \ref{ass}.  Thus, the algorithm cannot be applied in situations where cross-link interactions --- which would violate seperability --- exist in the network.  The applicability of the PBCD algorithm further hinges on the concavity of the marginal value function of the lower-level problem, which ensures Problem \eqref{alg1-sub3} is convex and can be solved globally.  This requirement may not be met in other mixed network design problems. A case in point is the capacity expansion problem, another classic network design application.  {Because the objective function in the lower-level problem depends on link capacities through the link travel time function, concavity cannot be guaranteed in this problem except for highly simplistic forms of the travel time function. }  Finding ways to generalize the PBCD algorithm for tackling other mixed network design problems constitutes an important direction for future research.

\bibliographystyle{apalike}
\begin{small}
\bibliography{references}
\end{small}

\newpage
\begin{appendix}
 \section{Omitted proofs}\label{appen}

\subsection{Proof of Proposition \ref{prop}}

The results can follow from \citet[][]{Guo2024-joc}. For completeness, we provide a proof here.

(i) Since $t_a(v_a)$ is strictly increasing with respect to $v_a$ for all $a\in \cA$, it follows that the lower-level objective function $f(u,v)$ is strictly convex with respect to $v$. Thus, $\cS(u)$ is a singleton for all $u$. Furthermore, by the well-known Danskin's theorem, it follows that the gradient of $\cV(u)$ exists and can be given by
\[
\nabla \cV(u) = \nabla_u f(u,v)|_{v=\cS(u)} = \cS(u).
\]
We next show that $\cS(u)$ is continuous. Since $\Omega$ is compact and independent of $u$, it is easy to verify that $\cV(u)$ is continuous with respect to $u$. Thus for any $u^k\to u^*$, we have
\[
\cV(u^k) = f(\cS(u^k),u^k) \to \cV(u^*).
\]
This means that all accumulation points of $\{\cS(u^k)\}_{k=1}^\infty$ belong to $\cS(u^*)$. Recalling the singleton property of $\cS(u^*)$, it follows that $\lim\limits_{k\to\infty} \cS(u^k) = \cS(u^*)$, indicating that $\cS(u)$ is continuous. Therefore $\cV(u)$ is continuously differentiable. 

(ii) We first show the concavity of $\cV(u)$. Let $u^1,u^2\in \Re^{|{\cA}|}$, and $\alpha \in [0,1]$. It suffices to show that
\[
\cV(\alpha u^1 +(1-\alpha)u^2) \ge \alpha \cV(u^1) +(1-\alpha)\cV(u^2).
\]
By the definition of $\cV(u)$, we have
\begin{eqnarray*}
&& \cV(u^1) = \min_{v'\in \Omega} f(u^1,v') \le f(u^1,v), \quad \forall v\in \Omega,\\
&& \cV(u^2) = \min_{v'\in \Omega} f(u^2,v') \le f(u^2,v), \quad \forall v\in \Omega.
\end{eqnarray*}
Multiplying the above inequalities by $\alpha$ and $1-\alpha$ respectively, and then summing them gives the following result:  for all $v\in \Omega$, 
\begin{eqnarray*}
\alpha \cV(u^1) +(1-\alpha)\cV(u^2) &\le& \alpha f(u^1,v) +(1-\alpha)f(u^2,v)\\
&=& \alpha \left[\sum\limits_{a\in \cA} \int_0^{v_a} t_a(x)dx + \sum\limits_{a\in {\cA}}u^1_a v_a\right] +(1-\alpha)\left[\sum\limits_{a\in \cA} \int_0^{v_a} t_a(x)dx + \sum\limits_{a\in {\cA}}u^2_a v_a\right]\\
&=& \sum\limits_{a\in \cA} \int_0^{v_a} t_a(x)dx + \sum\limits_{a\in {\cA}}(\alpha u^1_a +(1-\alpha)u^2_a) v_a \\
&=& f(\alpha u^1 +(1-\alpha)u^2,v),
\end{eqnarray*}
which indicates that
\[
\alpha \cV(u^1) +(1-\alpha)\cV(u^2) \le \min_{v\in \Omega}f(\alpha u^1 +(1-\alpha)u^2,v) = \cV(\alpha u^1 +(1-\alpha)u^2).
\]
Therefore $\cV(u)$ is concave with respect to $u$. This property, together with Assumption \ref{ass}, immediately implies the rest of the proof. $\qquad\blacksquare$
 
 \subsection{Proof of Proposition \ref{global}}

First from the definition of $u^*$, it follows that $\|u^*\|_0\le \kappa$, indicating that $u^*\in U_\kappa$ is feasible to Problem \eqref{projp}. Thus, by the definition of optimality to Problem \eqref{projp}, it suffices to show that
\[
\|u^*-z\|^2 \le \|u-z\|^2, \quad \forall u\in U_\kappa.
\]
Let $J_u = \{i\in \{1,\ldots,m\}: u_i = 0\}$ for each $u\in U_\kappa$. It is clear that for all $u\in U_\kappa$, the number of elements in $J_u$ satisfies $|J_u|\ge m-\kappa$ and hence
\begin{eqnarray*}
\|u-z\|^2 &=& \sum_{i\in J_u} (z_i)^2 + \sum_{i\notin J_u}(z_i-u_i)^2\\
&\ge& \sum_{i\in J_u} (z_i)^2 \ge   \sum_{i\in I_c} (z_i)^2 =\|u^{*}-z\|^2,
\end{eqnarray*}
where the last inequality follows from the fact that $I_c$ corresponds to the smallest $m-\kappa$ values of $\{|z_i|:i=1,\ldots,m\}$, and the last equality follows from the definition of $u^*$. $\qquad\blacksquare$

\subsection{Proof of Theorem \ref{sub-con1}}
We assume that for all $r\ge1$, $(u^r,z^r,v^r)$ is not a KKT stationary solution of the penalty approximation problem ${\rm (PA_\rho)}$. Otherwise, the algorithm will generate a finite sequence.

(i) We first show that $\Phi_\rho(u^{r+1},z^{r+1},v^{r+1})<\Phi_\rho(u^r,z^r,v^r)$ for all $r\ge 1$. By the optimality of $(u^{r+1},v^{r+1})$ to problem (\ref{alg1-sub}) and the optimality of $z^{r+1}$ to problem (\ref{alg1-sub3}), it follows that for all $r\ge1$,
\begin{eqnarray}
&& \Phi_\rho(u^r,z^r,v^r) \le  \Phi_\rho(u,z^r,v),\quad \forall u\in U_\kappa, v\in \Omega,\label{conv-opt1}\\
&& \Phi_\rho(u^r,z^{r+1},v^r) \le  \Phi_\rho(u^r,z,v^r),\quad \forall z\in U,\label{conv-opt2}\\
&& \Phi_\rho(u^{r+1},z^{r+1},v^{r+1}) \le \Phi_\rho(u,z^{r+1},v),\quad \forall u\in U_\kappa, v\in \Omega.\label{conv-opt3}
\end{eqnarray}
Since $u^r\in U_\kappa$, $z^r\in U$, and $v^r\in \Omega$, it follows that
\begin{equation}\label{eq:prof of theorem1}
\Phi_\rho(u^{r+1},z^{r+1},v^{r+1}) \le  \Phi_\rho(u^r,z^{r+1},v^r)\le \Phi_\rho(u^r,z^r,v^r),\ \forall r\ge 1.
\end{equation}
We next show that the strict inequality holds. To the contrary assume that $\Phi_\rho(u^{r+1},z^{r+1},v^{r+1})=\Phi_\rho(u^r,z^r,v^r)$. Then
\begin{equation}\label{equ-cont}
\Phi_\rho(u^{r+1},z^{r+1},v^{r+1}) = \Phi_\rho(u^r,z^{r+1},v^r)= \Phi_\rho(u^r,z^r,v^r).
\end{equation}
The first equality in \eqref{equ-cont}, together with the relation \eqref{conv-opt3}, implies that $(u^r,v^r)$ is a solution of Problem (\ref{alg1-sub}). Since Problem (\ref{d2}) is a strictly convex program, it follows that $v^{r+1}=v^r$. 
Similarly, the last equality in \eqref{equ-cont} and the relation \eqref{conv-opt2} imply that $z^r=z^{r+1}$ since Problem (\ref{alg1-sub3}) is a strongly convex program. Then by the choice rule of $u^r$ and $u^{r+1}$ as shown in Proposition \ref{global}, we know that $u^r=u^{r+1}$. Therefore, it follows 
$(u^r,z^r,v^r)=(u^{r+1},z^{r+1},v^{r+1})$. By \eqref{conv-opt1}, \eqref{conv-opt2}, and \eqref{conv-opt3}, it follows that
\begin{eqnarray*}
&& \Phi_\rho(u^r,z^r,v^r) \le \Phi_\rho(u,z^r,v),\quad \forall u\in U_\kappa, v\in \Omega,\\
&& \Phi_\rho(u^r,z^r,v^r) \le \Phi_\rho(u^r,z,v^r),\quad \forall z\in U.
\end{eqnarray*}
By Fermat's rule \citep[e.g.,][Theorem 10.1]{VA}, the above two inequalities imply that $(u^r,z^r,v^r)$ satisfies
\begin{eqnarray*}
0\in \nabla \Phi_\rho(u^r,z^r,v^r) + N_{U_\kappa}(u^r)\times N_{U}(z^r)\times N_{\Omega}(v^r).
\end{eqnarray*}
That is, it is a KKT stationary solution of ${\rm (PA_\rho)}$ which contradicts the assumption. Thus, the strict inequality in \eqref{monoto} holds.

(ii) Since both the sets $U$ and $\Omega$ are bounded, both $\{z^r\}_{r=1}^\infty$ and $\{v^r\}_{r=1}^\infty$ are  bounded. By the choice of $u^r$ as done in Proposition \ref{global}, the sequence $\{u^r\}_{r=1}^\infty$ is also bounded. Thus, the sequence $\{(u^r,z^r,v^r)\}_{r=1}^\infty$ is bounded. Noting that the objective function of problem \eqref{pp1p} is continuous, it follows that $\{\Phi_\rho(u^r,z^r,v^r)\}_{r=1}^\infty$ is bounded. This, together with the monotonicity in \eqref{monoto}, implies that the function value sequence must have a unique limit. Let $(u^*,z^*,v^*)$ be an accumulation point of $\{(u^r,z^r,v^r)\}_{r=1}^\infty$ and $T \subseteq \{1,2,\ldots\ldots\}$ be a subsequence such that $\lim_{r\in T\to \infty}(u^r,z^r,v^r)=(u^*,z^*,v^*)$.
Then by \eqref{eq:prof of theorem1}, we have
\begin{equation*}
\lim_{r\to\infty}\Phi_\rho(u^{r+1},z^{r+1},v^{r+1}) = \lim_{r\to\infty}\Phi_\rho(u^r,z^{r+1},v^r)= \lim_{r\to\infty}\Phi_\rho(u^r,z^r,v^r) = \Phi_\rho(u^*,z^*,v^*).
\end{equation*}
Using these relations and the continuity of $\Phi_\rho$, and taking limits on both sides of \eqref{conv-opt1} and \eqref{conv-opt2} respectively as $r\in T\to\infty$, we have
\begin{eqnarray}
&& \Phi_\rho(u^*,z^*,v^*) \le \Phi_\rho(u,z^*,v),\quad \forall u\in U_\kappa, v\in \Omega,\label{par2}\\
&& \Phi_\rho(u^*,z^*,v^*) \le  \Phi_\rho(u^*,z,v^*),\quad \forall  z\in U.\label{par1}
\end{eqnarray}
By Fermat's rule \citep[e.g.,][Theorem 10.1]{VA}, these two inequalities indicate that $(u^*,z^*,v^*)$ satisfies
\begin{eqnarray*}
0\in \nabla \Phi_\rho(u^*,z^*,v^*) + N_{U_\kappa}(u^*)\times N_{U}(z^*)\times N_{\Omega}(v^*),
\end{eqnarray*}
i.e., it is a KKT stationary solution of ${(\rm PA_{\rho})}$. $\qquad\blacksquare$

\subsection{Proof of Theorem \ref{mainthm1}}
Since $(u^k,z^k,v^k)$ is a globally optimal solution of ${(\rm PA_{\rho^{k}})}$, it follows that
\begin{equation*}
\Phi_{\rho^{k}} (u^k,z^k,v^k) \le \min_{v\in \cS(u),u\in U\cap U_\kappa}F(v) = F^*.
\end{equation*}
Thus by the explicit expression of $\Phi_{\rho^{k}} (u^k,z^k,v^k)$ and the definition of $F^l$, we have
\begin{equation}\label{opine}
F^l+\rho^{k}_1 (f(z^k,v^k)-\cV(z^k))+\rho^{k}_2 \|u^k-z^k\|^2 \le F(v^k)+\rho^{k}_1 (f(z^k,v^k)-\cV(z^k))+ \rho^{k}_2  \|u^k-z^k\|^2  \le F^*.
\end{equation}
Then we have
\begin{equation*}
f(z^k,v^k)-\cV(z^k)  \le \frac{F^*-F^l}{\rho^k_1},\ \|u^k-z^k\|^2 \le \frac{F^*-F^l}{\rho^k_2}.
\end{equation*}
When $\rho^k_1 = \gamma_1^{k-1}\rho_1^1 \ge\frac{F^*-F^l}{\varepsilon}$ and $\rho^k_2 = \gamma_2^{k-1}\rho_2^1 \ge\frac{F^*-F^l}{\varepsilon^2}$, i.e., the iteration number satisfies \eqref{number1}, it follows that
\begin{equation}\label{inequ12}
f(z^k,v^k)-\cV(z^k)\le \varepsilon,\ \|u^k-z^k\| \le \varepsilon.
\end{equation}
Furthermore, by the second inequality in \eqref{opine}, it follows that $F(v^k) \le F^*.$ This, together with the two inequalities in \eqref{inequ12}, implies the desired result. $\qquad\blacksquare$

\subsection{Proof of Theorem \ref{mainthm}}

By Theorem \ref{sub-con1}, it follows that for all $k\ge1$,
\begin{equation*}
\Phi_{\rho^k} (u^k,z^k,v^k) \le \min_{u\in U_\kappa,v\in \Omega} \Phi_{\rho^k}(u,z_0^k,v).
\end{equation*}
Then by the choice rule of the initial point for solving each penalty approximation problem, i.e., Step (3) of Algorithm \ref{mainalg1}, it follows that 
\begin{equation*}
\Phi_{\rho^k} (u^k,z^k,v^k) \le \min_{u\in U_\kappa,v\in \Omega} \Phi_{\rho^k}(u,z_0^k,v) \le \Upsilon.
\end{equation*}
Thus by the definition of $\Phi_{\rho}(u,z,v)$, it follows that
\begin{eqnarray*}
F^l+\rho_1^{k} (f(z^k,v^k)-\cV(z^k))+ \rho_2^k\|u^k-z^k\|^2 \le \Phi_{\rho^k} (u^k,z^k,v^k)  \le \Upsilon.
\end{eqnarray*}
Then we have
\begin{equation*}\label{inequ}
f(z^k,v^k)-\cV(z^k) \le \frac{\Upsilon-F^l}{\rho_1^{k}},\quad \|u^k-z^k\|^2 \le \frac{\Upsilon-F^l}{\rho_2^k}.
\end{equation*}
When $\rho_1^k = \gamma_1^{k-1}\rho^1_1 \ge\frac{\Upsilon-F^l}{\varepsilon}$ and $\rho_2^k = \gamma_2^{k-1}\rho^1_2 \ge\frac{\Upsilon-F^l}{\varepsilon^2}$, i.e., the iteration number satisfies \eqref{number}, it follows that
\begin{equation}\label{appfeas}
f(z^k,v^k)-\cV(z^k) \le \varepsilon,\ \|u^k-z^k\| \le \varepsilon.
\end{equation}
This means that $(u^k,z^k,v^k)$ is an approximately feasible solution of Problem \eqref{ppp1}. Then the stopping criteria in Step (2) of Algorithm \ref{mainalg1} are satisfied if the iteration number satisfies \eqref{number}. By the proof process of Theorem \ref{sub-con1} (i.e., \eqref{par2} and \eqref{par1}), the derived solution $(u^k,z^k,v^k)$ by Algorithm \ref{alg1} satisfies 
\begin{eqnarray}
&& \Phi_{\rho^k}(u^{k},z^k,v^k) \le  \Phi_{\rho^k}(u,z^k,v),\quad \forall u\in U_\kappa, v\in \Omega,\label{thproof-1}\\
&& \Phi_{\rho^k}(u^{k},z^{k},v^k) \le \Phi_{\rho^k}(u^k,z,v^k),\quad \forall z\in U.\label{thproof-2}
\end{eqnarray}
By the optimality of $(u^k,v^k)$ and $z^k$ as shown in \eqref{thproof-1} and \eqref{thproof-2}, and Fermat's rule \citep[e.g.,][Theorem 10.1]{VA}, it follows that
\begin{eqnarray*}
&& 0\in \nabla_z f(z^k,v^k) -\nabla \cV(z^k) + \mu^k + N_U(z^k),\\
&& 0\in -\mu^k + N_{U_\kappa}(u^k),\\
&& 0\in \nabla F(v^k) + \rho_1^k \nabla_v f(z^k,v^k) + N_\Omega(v^k).
\end{eqnarray*}
where $\mu^k = \frac{2\rho^k_2(z^k-u^k)}{\rho^k_1}$. These inclusions, together with the approximate feasibility \eqref{appfeas}, show that $(u^{k},z^k,v^k)$ is an approximate KKT solution of Problem \eqref{ppp1}. $\qquad\blacksquare$

\section{The path-based greedy algorithm for Problem (16)} \label{appen-b}

The path-based greedy algorithm proposed in \citet{xie2018greedy}  is a state-of-the-art method for solving the user equilibrium traffic assignment problem (UE-TAP). Since Problem \eqref{d2} has a convex objective function and shares the same constraints with UE-TAP, the path-based greedy algorithm can be easily adapted to solve it. To this end, we first write the explicit formulation of the objective function of Problem \eqref{d2} as follows:
\begin{align*}
 Z(v) &= F(v) + \rho_1 f(z^{r+1},v) \\ \nonumber
       &= \sum\limits_{a\in \cA} t_a(v_a) v_a + \rho_1 \sum\limits_{a\in \cA} \int_0^{v_a} (t_a(x)+z_a^{r+1})dx \\ \nonumber
         &= \sum\limits_{a\in \cA} \int_0^{v_a} (t_a(x)+xt_a'(x))dx + \rho_1 \sum\limits_{a\in \cA} \int_0^{v_a} (t_a(x)+z_a^{r+1})dx \\ \nonumber
         &= \sum\limits_{a\in \cA} \int_0^{v_a} \left( (1+\rho_1)t_a(x) + xt_a'(x) + \rho_1 z_a^{r+1} \right) dx. \nonumber
\end{align*}
For each link $a$, define 
\begin{align}
    c_a(v_a) &\equiv (1+\rho_1)t_a(v_a) + v_a t_a'(v_a) + \rho_1 z_a^{r+1}, \label{form:glc} \\
    c'_a(v_a) &\equiv (2+\rho_1)t'_a(v_a) + v_a t_a''(v_a),\nonumber
\end{align}
where $c_a(v_a)$ and $c'_a(v_a) $ can be viewed as a generalized link cost and its derivative respectively. Accordingly, Problem \eqref{d2} can be reformulated as 
\begin{equation*}\label{gtap}
\begin{array}{rl}
\rmin\limits_{v} & \hat{f}(v)=\sum\limits_{a\in \cA} \int_0^{v_a} c_a(x)dx\\[8pt]
\rst & v\in \Omega.
\end{array}
\end{equation*}
The above transformation reduces Problem \eqref{d2} to a form that resembles the path-based formulation for the traffic assignment problem. This enables us to solve the problem using the path-based greedy algorithm by simply replacing the original link cost function with the generalized link cost function \eqref{form:glc}. For details on the greedy algorithm, please refer to \citet{xie2018greedy}, and for its implementation, visit this \href{https://github.com/junxie016/Open-TNM/tree/main/Greedy}{GitHub repository}.

\section{The projected gradient algorithm for Problem (17)} \label{appen-c}

Problem \eqref{form:updatez} is reformulated as the following problem by incorporating the explicit formulation of $f(z,v^r)$:
\begin{eqnarray*}
\min\ g(z)=\rho_1\left(\sum\limits_{a\in \cA} \left( \int_0^{v_a^r}t_a(x)dx  + z_a v_a^r \right) - \cV(z)\right) + \rho_2\|z-u^r\|^2 \quad {\rm s.t.}\ z\in U,
\end{eqnarray*}
where $g:\Re^{|\cA|} \to \Re$ is a continuously differentiable convex function and $U$ is a convex set in $\Re^{|\cA|}$. Then the gradient of $g(z)$ can be computed by
\begin{equation*}
     \nabla g(z) = \rho_1(v^r - \cS(z))  + 2\rho_2(z - u^r),
\end{equation*}
where $\cS(z)$ is the gradient of $\cV(z)$ as given in Proposition \ref{prop} and represents the optimal solution of Problem \eqref{lp1} with $z$ in place of $u$. In light of the above results,
the projected gradient algorithm with the initial Barzilai–Borwein step size for Problem \eqref{form:updatez} is given as follows.

{
\setlist[itemize]{leftmargin=4.5em}
\setlist[itemize, 2]{leftmargin=2.5em}

\begin{alg}[The projected gradient algorithm]\label{pg}
Choose $0<\alpha_{\rm min}<\alpha_{\rm max}$, $\eta,\sigma\in (0,1)$, and an initial point $z^0\in U$. Set $\alpha_0=1$ and $k=0$.
\begin{itemize}
\item[Step (1):] If $\|P_U(z^k-\nabla g(z^k))-z^k\| =0$, then stop. Otherwise, go to Step (2).
\item[Step (2):] Set $\tau = \alpha_k$.
\begin{itemize}
\item[Step (2.1):] Set $z^+ = P_U(z^k-\tau \nabla g(z^k))$.

\item[Step (2.2):] If
    \begin{equation*}
    g(z^+) \le g(z^k) +\sigma \langle \nabla g(z^k), z^+-z^k\rangle,
    \end{equation*}
set $z^{k+1} = z^+$, $s^k = z^{k+1}-z^k$, $y^k = \nabla g(z^{k+1})-\nabla g(z^k)$ and go to Step (3). Otherwise, set $\tau = \eta \tau$ and go to Step (2.1).
\end{itemize}
\item[Step (3):] If $\langle s^k, y^k\rangle\le0$, then set $\alpha_{k+1} = \alpha_{\rm max}$. Otherwise, set
    \[
    \alpha_{k+1} = \min\big\{\alpha_{\rm max}, \max\big\{\alpha_{\rm min}, \frac{\langle s_k, s_k\rangle}{\langle s^k, y^k\rangle}\big\}\big\}.
    \]
    Let $k=k+1$ and go to Step (1).
\end{itemize}
\end{alg}
}

When applying Algorithm \ref{pg} for solving Problem \eqref{form:updatez} in this paper, we set $\alpha_{\rm min}= 10^{-20},$ $\alpha_{\rm max}=10^{20}$, $\eta=0.1$, $\sigma=0.01$, and $z^0 = 0$. Moreover, we terminate the algorithm once the iteration solution satisfies $\|P_U(z^k-\nabla g(z^k))-z^k\|\le 10^{-3}$.

\end{appendix}

\end{document}